\documentclass{article}
\usepackage{graphicx, subfigure, amsfonts, amsmath, amsthm}
\usepackage{a4wide, color}
\usepackage[top=1in, bottom=1.25in, left=1.25in, right=1.25in]{geometry}
\newcommand{\paf}[2]{\frac{\partial #1}{\partial #2}}

\newcommand{\discretized}{discretized}
\newcommand{\discretization}{discretization}
\newcommand{\discretizations}{discretizations}
\newcommand{\adjointness}{adjointness}
\newcommand{\adjoint}{adjoint}
\newcommand{\adjoints}{adjoints}
\newcommand{\Adjoints}{Adjoints}
\newcommand{\Nonnegative}{Nonnegative}
\newcommand{\nonnegative}{nonnegative}

\newcommand{\diag}{\mbox{diag}}

\newtheorem{thm}{Theorem}
\newtheorem{lem}[thm]{Lemma}

\begin{document}
\title{Symmetry-preserving discretizations of arbitrary order on structured curvilinear grids}
\author{Bas van 't Hof\footnote{bas.vanthof@vortech.nl. VORtech, P.O. Box 260, 2600AG Delft, The Netherlands.} \qquad Mathea~J. Vuik\footnote{Corresponding author. thea.vuik@vortech.nl. Telephone: +31 (0)152850125. VORtech, P.O. Box 260, 2600AG Delft, The Netherlands.}}

\maketitle

\begin{abstract}
Mathematical descriptions of flow phenomena usually come in the form of partial differential
equations. The differential operators used in these equations may have properties such as symmetry,
skew-symmetry, positive or negative (definite)-ness. 
Symmetry-preserving methods are such that the discretized form of the continuous differential operator exhibits the same properties as the continuous operator itself.  The use of symmetry-preserving \discretizations\ makes it possible to construct discrete models which allow all the 
manipulations needed to prove stability and (discrete) conservation properties in the same way they
were proven in the original continuous model.  Furthermore, these methods allow a \discretization\ of the 
continuous \adjoint\ which is at the same time the discrete \adjoint\ of the discrete forward model. Such
\adjoint\ models are not harder to code than the discrete forward model.

This paper presents a new symmetry-preserving discretization of arbitrary order on curvilinear structured grids. The key idea is to use mutually-adjoint sampling and interpolation operators to switch between the continuous and discrete operator. The novelty of this work is that it combines three important requirements for discretizations: first, the symmetry-preserving discretization is made for arbitrary order of accuracy; second, the method works for every structured curvilinear mesh; and third, the method can be applied to every continuous operator. This paper is the first in a series of papers that gradually extends the theory to a general approach.
\end{abstract}

\textbf{Symmetry-preserving discretizations, Energy conservation, Mass conservation, Cur\-vilinear grid}

\section{Introduction and motivation}
   Mathematical descriptions of flow phenomena usually come in the form of partial differential
equations. The differential operators used in these equations may have properties such as symmetry,
skew-symmetry, positive or negative (definite)-ness \cite{Ver03V}.  Proofs of stability and/or conservation
properties, such as conservation of mass, momentum, and energy, can often be constructed using the
symmetry and/or positiveness of the operators in the flow model \cite{Bla96SQ}. 

Computer simulations of a flow phenomenon require the \discretization\ of the flow properties, reducing
the number of values needed to represent the flow state from infinite to some large finite number.
In the resulting discrete model of the flow phenomenon, the differential operators have been
replaced by difference operators.  Unfortunately, not all properties of the differential operators
are automatically inherited by their discrete approximations.  The chain and product rules 
needed in the manipulation of nonlinear equations, for example, do not always work in discrete cases.
Moreover, symmetry and positiveness may be lost in the \discretization\ process, mass, momentum, and energy may not be conserved, aliasing errors can occur, and  duality and self-adjointness of the differential operators may be violated \cite{Bla96SQ, Kok09, Lip14MS}. 

The above difficulties play an important role when \adjoint\ models are considered.
Every (partial) differential model has an \adjoint\ model, as does every (discrete) difference
model. The \adjoint\ model is needed in optimization problems and is often used to find solutions for seismic models \cite{Bra13, Tro05TL}. A big dilemma in optimization is often 
the choice between the (discrete) \adjoint\ of the discrete forward model and the 
\discretization\ of the (continuous) \adjoint\ of the continuous forward model. \Adjoints\ of discrete
forward models often lead to very complicated and inefficient computer code \cite{Bra13}.

Symmetry-preserving methods are such that the discretized form of the continuous operator exhibits the same properties as the continuous operator itself \cite{Kok09}.  The use of symmetry-preserving \discretizations\ makes it possible to construct discrete models which allow all the 
manipulations needed to prove stability and (discrete) conservation properties in the same way they
were proven in the original continuous model.  Furthermore, these methods allow a \discretization\ of the 
continuous \adjoint\ which is at the same time the discrete \adjoint\ of the discrete forward model. Such
\adjoint\ models are not harder to code than the discrete forward model.

There are a variety of symmetry-preserving \discretizations\ available in the literature. In \cite{Hof12V}, an exhaustive overview is given of different techniques to obtain mass- or energy-conserving methods.
Typically, symmetry properties of differential operators are only automatically preserved in central-difference
approximations on uniform, rectilinear grids \cite{Kok09}. Although finite-volume methods can be used to construct
conservative discretizations for mass and momentum, it is in general not possible to also obtain energy conservation \cite{Hof12V}. 

 In \cite{Ver98V, Ver03V}, a fourth-order symmetry-preserving finite-volume method is constructed using Rich\-ardson extrapolation of a second-order symmetry-preserving method \cite{Vel92R}. The extension to collocated unstructured meshes is presented in \cite{Tri14LOPSV}, and an application can be seen in \cite{Ver06V}.  The extension to upwind discretizations was made in \cite{Vel08L}, and a discretization for the convective operator was found in \cite{Kok06}. In \cite{Kok09}, the method is extended to non-uniform curvilinear structured grids by deriving a discrete product rule. Furthermore, a symmetry-preserving method that conserves mass and energy for compressible flow equations with a state equation is described in \cite{Hof12V}. For rectilinear grids, this method works well, but it is challenging to let this method work for unstructured grids. 
 
 Another option to preserve symmetry is to use discrete filters to regularize the convective terms of the equation \cite{Tri11V, Leh12BRPSO}. The combination of a symmetry-preserving discretization and regularization for compressible flows is studied in \cite{Roz15VKV}. 
 
 Mimetic finite-difference methods also mimic the important properties of differential operators. An interesting review is given in \cite{Lip14MS}, and recently, a second-order mimetic discretization of the Navier-Stokes equations conserving mass, momentum, and kinetic energy was presented in \cite{Oud16HVH}. The mimetic finite-difference method uses algebraic topology to design and analyze compatible discrete operators corresponding to a continuous formulation \cite{Boc06H, Kre11PG}. In order to construct a discrete de Rham complex, certain conditions on reconstruction and reduction operators are imposed: they should be conforming, which means that the reconstruction is a right inverse of the reduction \cite{Boc06H}, they should be constant preserving \cite{Bre10B}, and the interpolation operator should commute with the differential operator \cite{Bre10B}. In \cite{Bre14BM} a nice overview of mimetic methods is given. Discrete exterior calculus (DEC) is also related to these mimetic approaches \cite{Hir03}.

 The above papers all have their own advantages and disadvantages. In general, they are not applicable for all orders, operators, and meshes. The current paper presents a new discretization method that can handle these requirements simultaneously.
 
In this work, we present a new symmetry-preserving discretization of arbitrary order on curvilinear structured grids. The key idea is to use mutually-adjoint sampling and interpolation operators to switch between the continuous and discrete operator. The novelty of this work is that it combines three important requirements for discretizations: first, the symmetry-preserving discretization is made for arbitrary order of accuracy; second, the method works for every structured curvilinear mesh; and third, the method can be applied to every continuous operator. This paper is the first in a series of papers that gradually extends the theory to a general approach.

The outline of this paper is as follows: in Section \ref{sec:background}, we present the relevant background information on positivity and symmetry, and in Section \ref{sec:symmpres}, we introduce our new symmetry-preserving discretization. The effectivity of this new method is presented in Section \ref{Sec: example} for the wave equation. We conclude with a discussion of our method and future work in Section \ref{sec:conclusion}.

\section{Background}\label{sec:background}
Before presenting our new symmetry-preserving \discretization, a short discussion is needed about the definitions we use for positiveness and 
symmetry preservation.

\subsection{Inner products and positivity}
In the linear space of continuous, square-integrable functions in a domain $D$, 
the standard scalar product is 
\[
 \langle f,g\rangle = \int_D f(\vec x)^Tg(\vec x)~\mbox dD.
\]

In finite-dimensional spaces, any scalar product has to be given by
\begin{eqnarray}
\langle x,y\rangle = x^T {\bf Q} y
\label{eq: discrete scalar product}
\end{eqnarray}
for a matrix ${\bf Q}$ which is symmetric
and positive definite with respect to the standard scalar product. 

Let us define operator ${\cal A}$, having domain $V$ and co-domain $W$: if $x \in V$, then the {\em image} ${\cal A}x \in W$. Operators for which the domain and co\-domain are the same ($W=V$) are called \emph{square}: they
images 'live' in the same space as the arguments. Only square operators can be positive or negative (definite).

We distinguish the following categories:

\mbox{}\hfill
\begin{tabular}{llll}
{\bf \nonnegative}: && $\langle x,{\cal A}x\rangle \geq 0$ & for all $x$ in $V$,\\
{\bf positive}:    && $\langle x,{\cal A}x\rangle >    0$ & for all $x\neq 0$ in $V$,\\
{\bf positive definite}: & $\epsilon>0$ exists so that &$\langle x,{\cal A}x\rangle \geq \epsilon \langle x,x\rangle \geq 0$ & for all $x$ in $V$,\\
\end{tabular}
\hfill\mbox{}\\
and similarly, nonpositive and negative (definite) operators are defined. 
Note that 
all positive definite operators are positive, and
all positive operators are \nonnegative.

In finite-dimensional spaces, every positive operator is also positive definite. \Nonnegative\
operators which are not positive must be singular: they must have a non-empty null space.

\subsection{Adjoint operators and symmetry}
The words \adjoint\ and transpose are often used interchangeably. In many cases, 
the difference between the transpose and the \adjoint\ is not important, but in the case of 
symmetry preservation, we have to specify the notion of \adjoint\ operators a little more precisely.

The definition of the \adjoint\ always depends on the definitions of the inner product
$\langle\cdot,\cdot\rangle_V$ for the domain, and the inner product $\langle\cdot,\cdot\rangle_{W}$ for the co\-domain of the operator
${\cal A}$. 
The {\em \adjoint}\ ${\cal A}^*$ of the operator ${\cal A}$ is the unique operator for which 
\begin{equation}
   \langle y,{\cal A}x\rangle_{W} =\langle{\cal A}^* y, x\rangle_{V} \quad \mbox{for all $x$ in $V$ and all $y$ in $W$}. \label{eqn:adjoint}
\end{equation}

In finite-dimensional, real spaces $V=\mathbb{R}^n$, when using the standard scalar product $\langle x,y\rangle=x^Ty$,
the \adjoint\ of a matrix is its transpose: $A^*=A^T$.

If the forward and \adjoint\ 
operators are the same (${\cal A}^* = {\cal A}$), we call the operator ${\cal A}$ 
\emph{symmetric}; if they are each other's opposites (${\cal A}^*=-{\cal A}$), the operator 
is called \emph{skew-symmetric}.

\section{Symmetry-preserving \discretizations}\label{sec:symmpres}
In this section, we propose a simple strategy for the construction of symmetry-preserving calculations for any
differential operator ${\cal A}$. Let ${\cal J}_p$ be an interpolation operator that maps from the discrete field to the continuous field, and ${\cal S}_p$ be the sampling operator that produces discrete values from a continuous function \cite{Vet14KG}. 
We shall write the interpolated fields using italic letters and sampled fields with bold-faced
letters, for example, $f := {\cal J}_p~{\bf f}$, and ${\bf g} = {\cal S}_p g$.

The continuous operator ${\cal A}$ is applied to the continuous
field obtained from interpolation of the discrete field using ${\cal J}_p$, 
and the result is mapped back using ${\cal S}_p$, which leads to the discrete
operator ${\bf A} := {\cal S}_p~{\cal A}~{\cal J}_p$.

\textbf{Definition 1.} The sampling operator and the interpolation operator are called {\em mutually adjoint} if ${\cal S}_p = {\cal J}^*_p$.

Mutual \adjointness\ of the sampling and interpolation operator
means that when the inner product of a continuous field $f$ and a sampled field ${\bf g}$ 
is calculated, it does not matter whether the sampled field is interpolated and the 
continuous inner product is used, or whether the continuous field is sampled and the discrete inner
product is used: the result is the same.

\begin{lem}
Let ${\bf f} = {\cal J}_p f$ and ${\bf g} = {\cal J}_p g$. If the sampling operator and the interpolation operator are mutually adjoint,
then the inner products in the discrete and continuous spaces are the same:
\[
\langle {\bf A}^*{\bf f}, {\bf g}\rangle_{\bf p} = 
\langle {\bf f}, {\bf A} {\bf g}\rangle_{\bf p} = 
\langle f,{\cal A}~g\rangle_{p} = 
\langle {\cal A}^* f,g\rangle_{p},
\]
and the symmetry properties of ${\cal A}$ will be preserved in the \discretization\ ${\bf A}$.

Moreover, when $\langle f,f\rangle_{p} \geq \epsilon \langle {\bf f},{\bf f}\rangle_{\bf p}$, then
${\bf A}$ will also inherit all positiveness properties from ${\cal A}$.
\end{lem}
\textbf{Proof.} When we apply the mutual adjointness of ${\cal S}_p$ and ${\cal J}_p$ and definition (\ref{eqn:adjoint}), we find
 \begin{align*}
  \langle {\bf f}, {\bf Ag}\rangle_{\bf p} = \langle {\bf A}^*{\bf f}, {\bf g}\rangle_{\bf p} \hspace{-0.1mm}=\hspace{-0.1mm} \langle ({\cal S}_p {\cal AJ}_p)^*{\bf f}, {\bf g}\rangle_{\bf p}\hspace{-0.1mm} &\hspace{-0.1mm}=\hspace{-0.1mm} \langle ({\cal J}_p^* {\cal AJ}_p)^*{\bf f}, {\bf g}\rangle_{\bf p} \\
  &= \hspace{-0.1mm}\langle {\bf f}, {\cal J}_p^* {\cal AJ}_p{\bf g}\rangle_{\bf p}\hspace{-0.1mm} =\hspace{-0.1mm} \langle {\cal J}_p{\bf f}, {\cal AJ}_p{\bf g}\rangle_{\bf p} \hspace{-0.1mm}= \hspace{-0.1mm}\langle f,{\cal A}~g\rangle_{p} \hspace{-0.1mm} = \hspace{-0.1mm} \langle {\cal A}^* f,g\rangle_{p}.
 \end{align*}
For the second statement, we assume that ${\cal A}$ is positive definite. Then 
\[
 \langle {\bf f}, {\bf Af}\rangle_{\bf p} = \langle {\bf f}, {\cal S}_p {\cal A} {\cal J}_p {\bf f}\rangle_{\bf p} = \langle {\cal J}_p {\bf f}, {\cal AJ}_p {\bf f}\rangle_p = \langle f, {\cal A} f\rangle_p \geq \epsilon_1\langle f, f\rangle_p \geq \epsilon_2 \langle {\bf f}, {\bf f}\rangle_{\bf p}:
\]
$\bf{A}$ is positive definite as well.
\qed

To see what the mutual \adjointness\ actually means for the sampling and interpolation operators,
they are both written in a more explicit form \cite{Vet14KG}:
\begin{equation}
({\cal J}_p{\bf g})(\vec x) = \sum_i {\bf g}_i w_i(\vec x) \quad , \quad
({\cal S}_pf)_i = \int_V f(\vec x)  s_i(\vec x) ~\mbox dV. \label{eqn:JS}
\end{equation}
Here, $w_i$ are the interpolation base functions, and $s_i$ are the sampling functions.
The integral of the sampling function should be 1, because that means that the sampling of a constant field is exact. 

\begin{lem}
Mutual adjointness of the interpolation and sampling operators is equivalent to 
the following definition of the interpolation function: 
\begin{equation}
     w_i(\vec x) = \sum_j s_j(\vec x) {\bf Q}_{ij},
    \label{eq: Lemma adjointness}
\end{equation}
which means that interpolation functions can be computed if the sampling functions are chosen.

In equation (\ref{eq: Lemma adjointness}), ${\bf Q}$ is the matrix belonging to the discrete inner 
product (equation (\ref{eq: discrete scalar product})).
\end{lem}
\textbf{Proof.} Mutual adjointness of the interpolation and sampling operator means that 
\[
 \langle f,{\cal J}_p{\bf g}\rangle_{p} = \left\langle {\cal S}_p f,{\bf g}\right\rangle_{\bf p}.
\]
Here,
 \[
  \langle f,{\cal J}_p{\bf g}\rangle_{p} = \int_V \sum_i w_i(\vec x) {\bf g}_i f(\vec x) dV,
 \]
and, since ${\bf Q}$ is symmetric, we find
\[
 \left\langle {\cal S}_p f,{\bf g}\right\rangle_{\bf p} \hspace{-0.1mm} =\hspace{-0.1mm} \sum_j \sum_i ({\cal S}_p f)_j {\bf Q}_{ji} {\bf g}_i \hspace{-0.1mm}= \hspace{-0.1mm}\sum_i \sum_j \left( \int_V s_j(\vec x) f(\vec x) dV \right) {\bf Q}_{ij} {\bf g}_i \hspace{-0.1mm}=\hspace{-0.1mm} \int_V \sum_i {\bf g}_i \sum_j  s_j(\vec x) {\bf Q}_{ij} f(\vec x) dV,
\]
such that the condition $\langle f,{\cal J}_p{\bf g}\rangle_{p} 
= 
\left\langle {\cal S}_p f,{\bf g}\right\rangle_{\bf p} $ results in  $w_i(\vec x) = \sum_j s_j(\vec x) {\bf Q}_{ij}$.

\qed

If the integration matrix {\bf Q} had off-diagonal elements, this would mean that the interpolation
base functions $w$ become the linear combination of multiple sampling functions. Since we wish both the 
sampling and interpolation base functions to be zero except in a small region near a grid point, it makes
no sense that the interpolation base function should be nonzero in a larger region than the sampling 
function. Therefore, we shall expect the integration matrix ${\bf Q}$ to be a 
diagonal matrix. From here on, we will, therefore, use the notation $\diag({\bf Q})$ for the integration matrix, 
and the matrix entries ${\bf Q}_i$ need only one row/column index.
The interpolation base functions become scaled versions of the sample functions:
\[
    w_i(\vec x) = {\bf Q}_i s_i(\vec x).
\]

Since $\int_V s_j(\vec x) dV = 1$, we find that ${\bf Q}_i = \int_{V_i} w_i(\vec x) dV$.
Now that ${\bf Q}$ is diagonal, we can easily choose the interpolation functions and compute the corresponding sampling functions. 

\subsection{Interpolation functions on a uniform 1D grid}
The simplest set of interpolation base functions is the case of an infinite, uniform one-dimensional grid with unit grid
distance $\Delta x=1$. In such a case, only one base function is enough to construct all the other
ones, because they are found by translation:  $w_j(x) = w_0(x-j)$.  
The rest of this Section describes a method for the construction of the 
interpolation spline $w_0$. 

To obtain a unique spline, the following parameters are chosen:
\begin{description}
   \item{\tt nSpan}: the span of the function's support: $w_j(x)=0$ for all $x < j-{\tt nSpan}$ and for all $x>j+{\tt nSpan}$;
   \item{\tt nCont}: the number of continuous derivatives of the interpolation spline (internally), and the number of zero derivatives of the spline at the boundaries;
   \item{\tt Order}: the spline is a piecewise polynomial of order {\tt Order};

        This parameter is not very important: for sufficiently large orders, the interpolation base function 
        no longer depends on it.

   \item{\tt nConsist}: interpolation of all polynomials up to order {\tt nConsist} will be exact, and the interpolation will converge with order {\tt nConsist};

   \item{\tt wmax}: largest grid wavenumber for which the interpolation of the function $f(x) =
\exp(i x~\tt wmax)$ is accurate.
\end{description}
The combination of linear constraints ({\tt nCont}, {\tt nConsist})
and linear least squares equations ({\tt wmax}) leads to a unique set of interpolation base functions.  The interpolation obtained has a formal order of accuracy given by {\tt
nConsist}, and will be accurate for grid wavenumbers up to the given maximum. 

In Table \ref{tab:splines}, the parameter choices for three different interpolation splines are given. These splines are used for the computations in the rest of this paper. 

\begin{table}[ht!]
\centering
\caption{Parameters for three different interpolation splines that are used in the rest of this paper. The three splines all satisfy {\tt nSpan} $ = 3$, {\tt nCont} $=1$, {\tt Order } $=11$. } \label{tab:splines}

\vspace{0.2cm}
\begin{tabular}{c|ccc}
& coarse & medium & fine \\
\hline
{\tt wmax} & 0.9 & 0.6 & 0.5 \\
{\tt nConsist} & 3 & 3 & 4
\end{tabular}
\end{table}

\begin{figure}[ht!]
\centering
\subfigure[Interpolation splines]{\includegraphics[height=6cm]{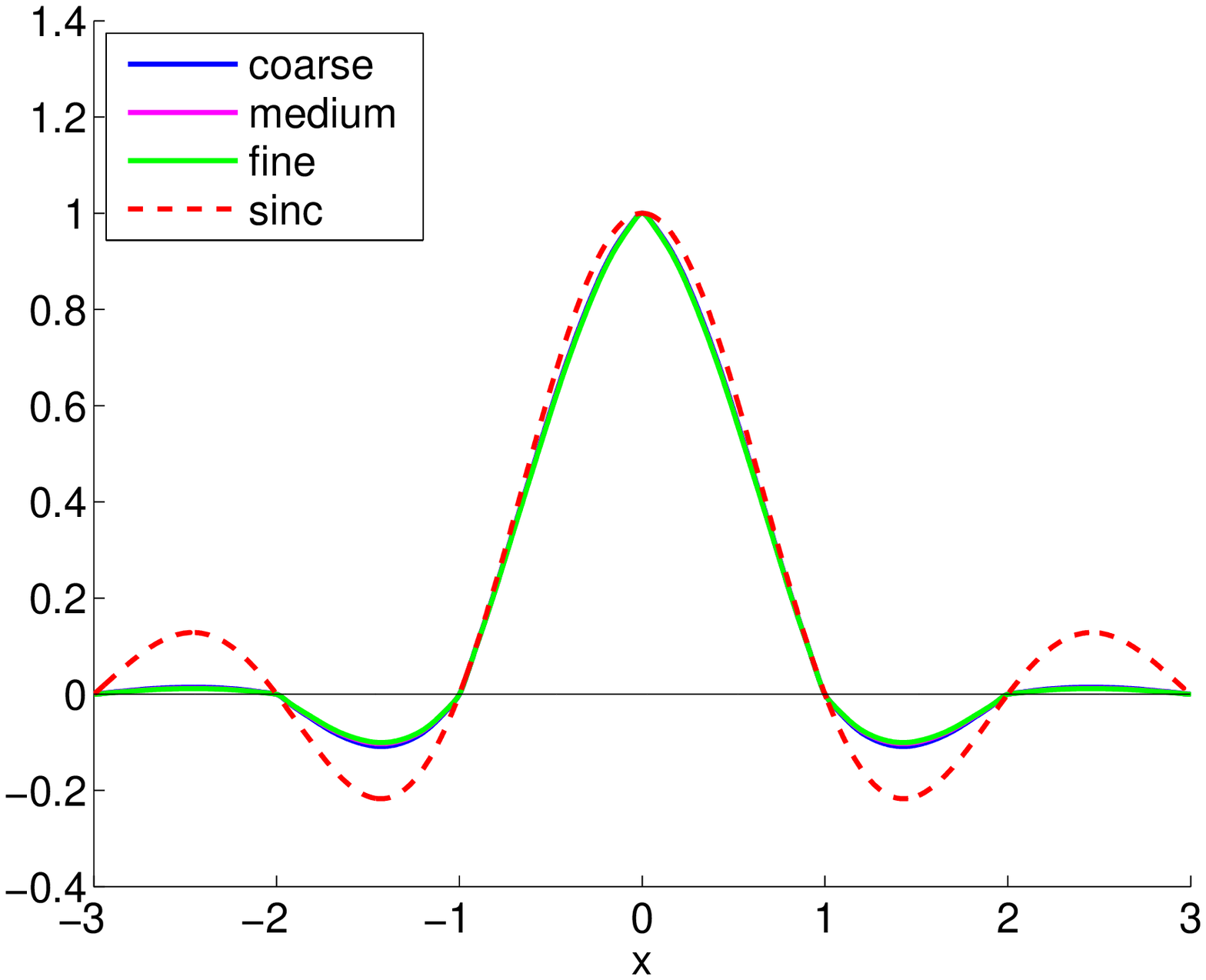}}
\subfigure[Largest error interpolating $e^{i\omega x}$]{\includegraphics[height=6cm]{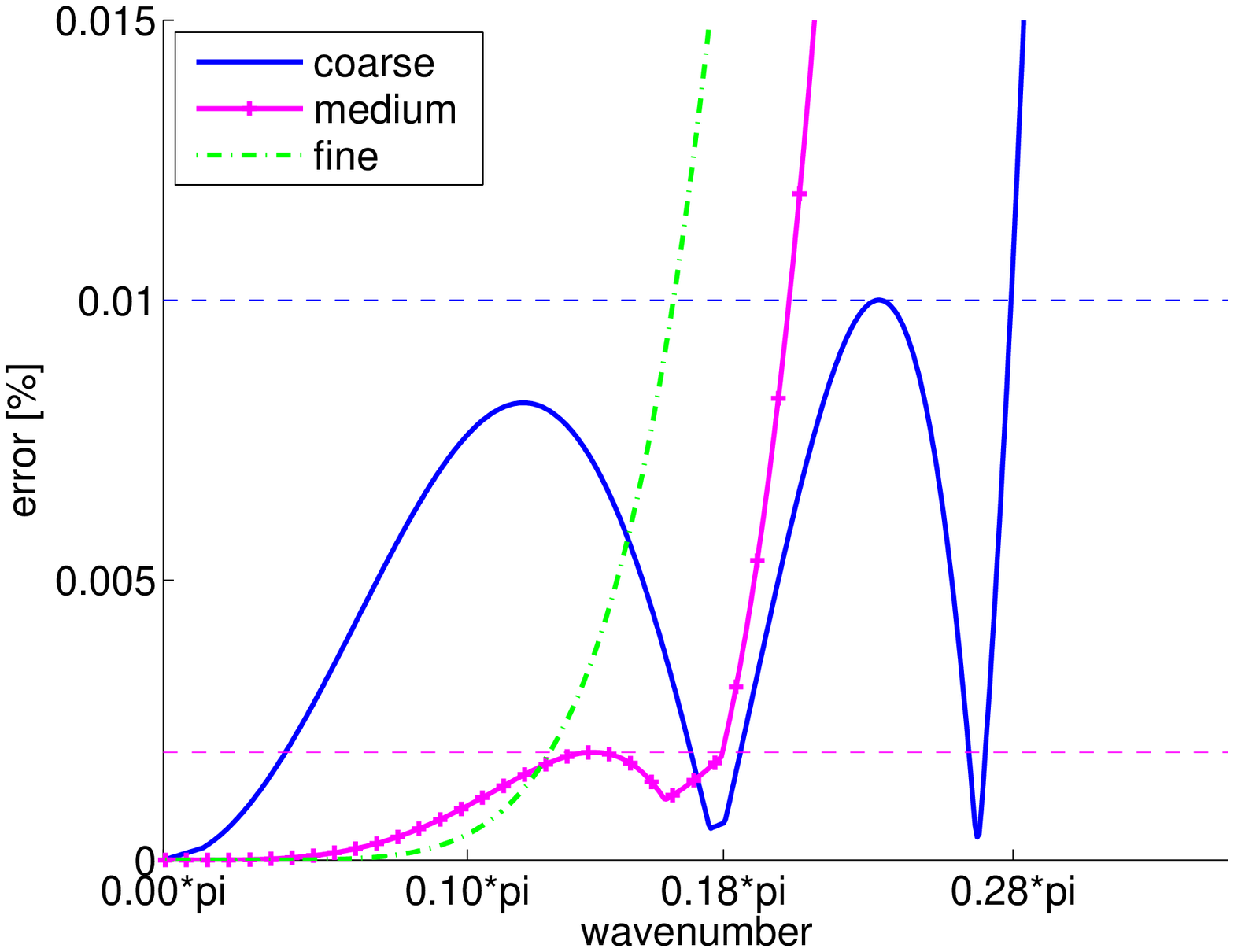}\label{fig:maxerror3}} \\
\caption{Interpolation base functions belonging to {\tt nSpan}=3, {\tt nCont}=1, {\tt Order}=11, for different values of {\tt wmax} and {\tt nConsist}.
The base function is compared to the {\tt sinc} function, 
which has all desired properties, but has infinite support. At the right, the largest interpolation errors over the domain are shown.}
\label{Fig: Wavelet9}
\end{figure}

Examples of interpolation base functions are shown in Figure \ref{Fig: Wavelet9}. 
The grid wavenumber $\omega$ corresponds to the number of grid points per wavelength in the grid \cite{Kim96L, Kuo14LL}. 
For a given problem, which has a given Fourier spectrum, refining the grid has the effect of reducing the grid wavenumber.  A sufficiently accurate solution will be found when the interpolation errors for all relevant grid wavenumbers are small.
Therefore, Figure \ref{fig:maxerror3} illustrates how the accuracy of the interpolations may be tuned for a specific problem. We discuss two possible 
scenarios:
\begin{itemize}
\item {\bf Modest accuracy: low convergence order}

   In the first example, it will be assumed that an accuracy of 0.01\% is required.  Figure \ref{fig:maxerror3} shows an interpolation of third order accuracy, that is sufficiently accurate for grid wavenumbers up to $0.28 \pi$.  The interpolation of fourth order accuracy requires refinement of the grid until grid wavenumbers are below $0.21\pi$.  With the lower-order interpolation, 36\% fewer grid points are needed in each direction.

\item {\bf Very accurate: high convergence order}

   For an accuracy of 0.002\%, the 4th order interpolation allows grid wavenumbers up to $0.18\pi$, while the lower-order interpolation needs refinement until grid wavenumbers are below $0.04\pi$. In this case, the higher-order interpolation requires 4.5 times fewer grid points in each direction than the lower-order interpolation. 

   The error in $\omega = 0$ always equals zero.
   The number of zero derivatives in $\omega = 0$ is equal to the convergence order of the interpolation spline, so for a very accurate solution, an interpolation is needed that has the largest possible order of convergence. 
\end{itemize}

\label{Sec: uniform 1D}
\subsection{Interpolation on structured, nonlinear 3D grids}
When using structured, nonlinear three-dimensional grids, we assume that some mapping exists 
between the $(x,y,z)$-locations of the physical domain and the $(\xi,\eta,\zeta)$-locations 
in a domain which we shall call {\em array space}. The physical domain points are found from the 
array-space locations by applying functions:
\[
 \left( \begin{array}{ll}
     x(\xi,\eta,\zeta) \\
     y(\xi,\eta,\zeta) \\
     z(\xi,\eta,\zeta) \end{array}\right) = \vec x(\xi,\eta,\zeta) = \vec x(\vec \xi).
\]
The function $\vec x(*)$ is assumed to be 'reasonably smooth'.
The grid points are found by entering integer values for the array-space coordinates $\xi$, $\eta$ and $\zeta$:
\[
\vec x_{ijk} := \vec x(i,j,k).
\]
The base function for uniform 1D grids calculated in Section
\ref{Sec: uniform 1D} is now combined into interpolation base functions and sampling functions 
according to
\begin{align*}
w_{ijk}(\vec x(\xi,\eta,\zeta)) &:= w_0(\xi-i) w_0(\eta-j) w_0(\zeta-k), \\
s_{ijk}(\vec x) &:= \frac{w_{ijk}(\vec x)}{\int_V w_{ijk}(\vec x)~\mbox dV}.
\end{align*}
This choice of the interpolation base functions secures the exact interpolation of constant fields, which is essential in the proof of discrete mass conservation in the example 
of Section \ref{Sec: example}.

The diagonal integration matrix ${\bf Q}$ is given by
\[
{\bf Q}_{ijk} := \int_V w_{ijk}(\vec x) ~\mbox dV.
\]

\section{Example: wave equation}
\label{Sec: example}
In this section, we show the effectivity of our new method by investigating the wave equation.
\subsection{The wave equation and its symmetry}
A simple illustration of the effects of symmetry preservation involves the 
wave equation:
\begin{eqnarray}
 \paf{{}^2 p}{t^2} = \nabla^2 p,  
\label{eq: wave eq}
\end{eqnarray}
where $p(\vec x, t)$ is the pressure.
Initial conditions for $p$ and $\partial p/\partial t$ as well as one boundary condition are needed to determine 
a solution of the wave equation.
Mass conservation requires initial conditions for $\partial p/\partial t$ that cause the initial time 
derivative of mass to be zero.
Therefore, we consider only initial conditions for the time derivative with a zero integral
\[
 \int_V \paf{{} p}{t}(\vec x,0) ~\mbox dV  = 0.
\]

The solution of the wave equation should conserve at least two quantities: the mass $M$ and the energy $E$, given by
\[
  M(t) := \int_V p (\vec x,t) ~\mbox dV \quad , \quad 
  E(t) := \frac 12\int_V \left( \left(\paf pt\right)^2 + |\nabla p|^2 \right) ~\mbox dV.
\]

Mass and energy are conserved in the sense that any changes over time can be expressed as 
the result of boundary terms called {\em fluxes}.
Using the divergence theorem, 
the change in energy is given in terms of {energy fluxes}:
\begin{align*}
   E'(t) &= 
        \int_V \left( \paf{{}^2 p}{t^2}\paf pt + \nabla p \cdot \nabla \paf pt \right)~\mbox dV
=
        \int_V \left(\nabla^2 p\paf pt + \nabla p \cdot \nabla \paf pt \right)~~\mbox dV
\nonumber \\ &=
        \int_V \nabla \cdot(\nabla p\paf pt) ~\mbox dV
=
        \oint_{\delta V} \left( \nabla p\paf pt ~\right) \cdot \vec n~\mbox dS.
\end{align*}

In formulation (\ref{eq: wave eq}) of the wave equation, it is not the first but 
the second-order time derivative  of the mass that
can be expressed in terms of {\em mass fluxes}:
\[
   M''(t) 
=
        \int_V \paf{{}^2 p}{t^2}~\mbox dV =  \int_V \nabla^2 p~\mbox dV
=
        \oint_{\delta V} \nabla p \cdot \vec n~\mbox dS.
\]
A consequence of these expressions of the time derivatives of mass and energy is that the energy and 
mass remain unchanged in unbounded or periodical domains.

Especially the case of energy conservation is interesting, because the expression used for the 
energy contains a norm of the solution. Hence, energy conservation means that the norm of the solution remains 
constant, which implies stability of the solution: energy conservation implies stability.

Conservation of mass and energy can also be proven using \adjoints\ and symmetry.
To do this, we write mass and energy as scalar products:
\[
   M(t) = \langle 1,p\rangle_{p} \quad,\quad
   E(t) = \frac 12\left\langle \paf pt,\paf pt\right\rangle_{p} + \frac 12 \left\langle \nabla p,\nabla p\right\rangle_{p}.
\]
In periodic domains, the divergence and the negative gradient are mutually adjoint ($(\nabla\cdot)^* = -\nabla$) \cite{Hym97S}, and so we get
\begin{align*}
   E'(t) &= \left\langle \paf pt,\nabla \cdot \nabla p\right\rangle_{p} + \left\langle \nabla p,\nabla \paf pt\right\rangle_{p} = 
   -\left\langle \nabla \paf pt,\nabla p\right\rangle_{p} + \left\langle \nabla p,\nabla \paf pt\right\rangle_{p} = 0,
    \\
   M''(t) &= \langle 1,\nabla^2 p\rangle_{p}=-\langle \nabla 1,\nabla p\rangle_{p} = 0.
\end{align*}
\subsection{Discrete model}
In the discrete case, we have precisely the same results. The discrete wave equation will be 
given by
\[
 {\bf p}''(t) = {\bf A} ~{\bf p} =: {\cal S}_p \nabla^2 {\cal J}_p {\bf p}.
\]
Discrete mass {\tt M} and energy ${\tt E}$ will be given by
\[
   {\tt M}(t) := \left\langle {\bf 1},{\bf p}\right\rangle_{\bf p}\quad,\quad
{\tt E}(t) := 
\frac 12
\left\langle  {\bf p}', {\bf p}'\right\rangle_{\bf p} - 
\frac 12 \left\langle  {\bf p}, {\bf A}~{\bf p}\right\rangle_{\bf p}.
\]
and their time derivatives are zero, using the same steps as in the continuous proof, and the mutual adjointness of the sampling and interpolation operator: 
\begin{align*}
   {\tt M}''(t) &:= \left\langle {\bf 1},{\bf A}~{\bf p}\right\rangle_{\bf p} = \left\langle {\bf 1},{\cal S}_p\nabla^2~{\cal J}_p~{\bf p}\right\rangle_{\bf p} =
 \left\langle {\cal J}_p {\bf 1},\nabla^2~{\cal J}_p~{\bf p}\right\rangle_{p} = 
-\left\langle \nabla 1,\nabla~{\cal J}_p~{\bf p}\right\rangle_{p} = 0.
 \\
{\tt E}'(t) &:= 
\left\langle  {\bf p}', {\bf A}~{\bf p}\right\rangle_{\bf p} - 
\left\langle  {\bf p}, {\bf A}~{\bf p}'\right\rangle_{\bf p} = 0.
\end{align*}
The proof of mass conservation requires the perfect interpolation of the constant field: ${\cal J}_p~{\bf 1} = 1$.

\subsection{Numerical results}

\begin{figure}[ht!]
\centering
\subfigure[Uniform mesh]{\includegraphics[scale = 0.4]{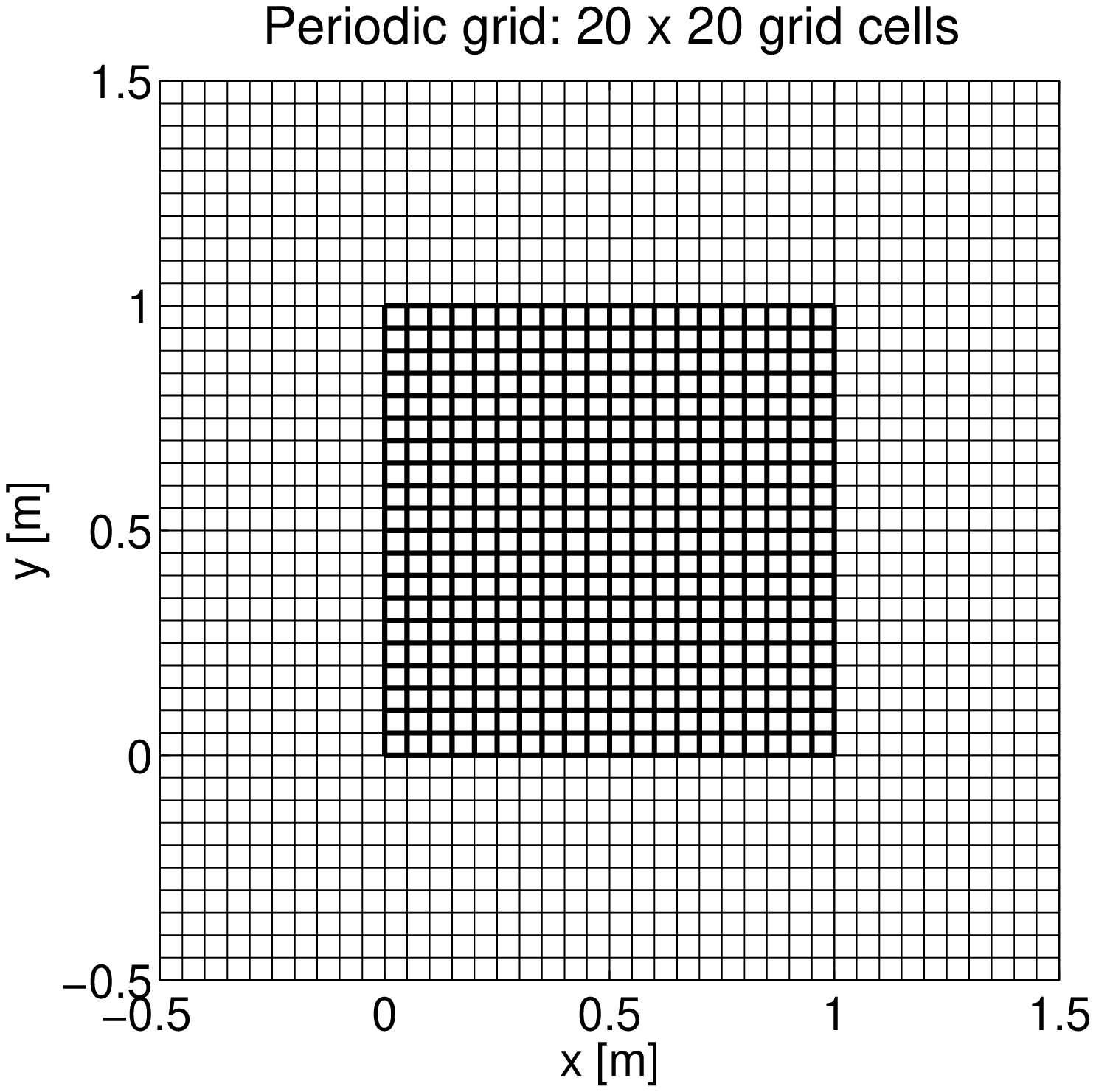}}
\subfigure[Curvilinear mesh]{\includegraphics[scale = 0.4]{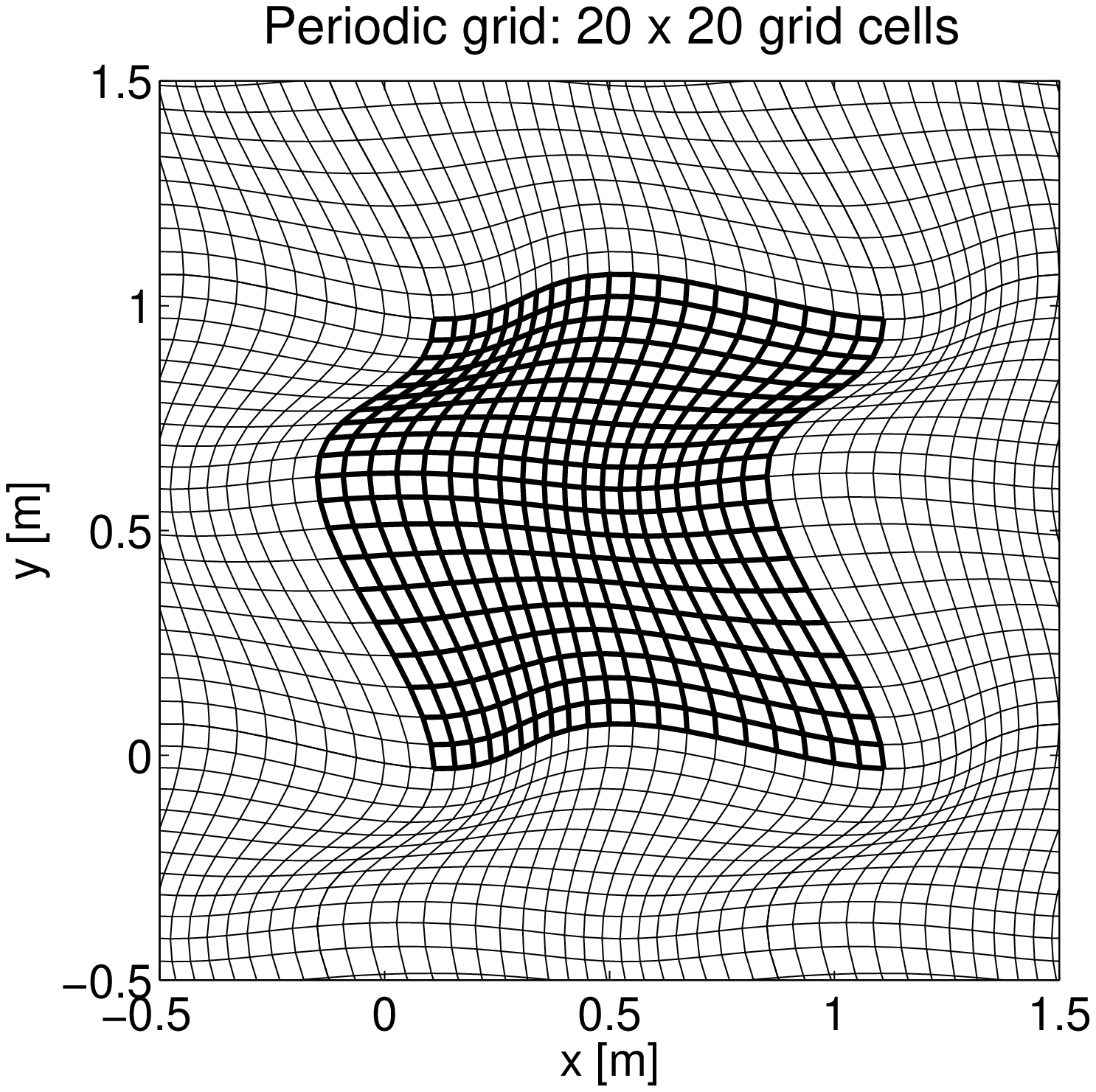}}
\caption{2D grid on $[0,1] \times [0,1]$ used for the \discretization\ of the wave equation. An impression of the periodicity is given by extending the mesh on each boundary. }
\label{Fig: grid}
\end{figure}

The wave equation is \discretized\ on a uniform and a 2D curvilinear grid for the unit square
with periodic boundaries, shown in Figure \ref{Fig: grid}. This results in the following equation: ${\bf p}^{\prime \prime} = {\bf A p}$, where ${\bf A}$ is given by ${\bf A} = {\cal S}_p \nabla^2 {\cal J}_p $, and definition (\ref{eqn:JS}) is applied.
Note that the discretization matrix ${\bf A}$ belonging to the wave equation is negative definite, and hence, it has negative eigenvalues. 

A Runge-Kutta time-integration method is used, where it is verified that the time integration is accurate enough such that it has no influence on the results. 
The initial conditions are chosen such that the exact solution is a one-dimensional, 
traveling Gaussian wave, given by
\[
   p_{exact}(x,y,t) = \exp\left( -\left(\frac{\mbox{mod}\left(x-y-t\sqrt 2+\frac 12,1\right)-\frac 12}{\sqrt{0.03}}
\right)^2 \right).
\]
The exact solution is sampled directly to obtain a reference solution:
\[
   \left({\bf p}_{ref}(t)\right)_{ij} := 
   p_{exact}(x_{ij},y_{ij},t).
\]
Initial conditions for {\bf p} and ${\bf p}'$ 
are taken from the reference solution, but a correction is applied 
to the time derivative ${\bf p}'$:
\[
    {\bf p}_0 := {\bf p}_{ref}(0) \quad,\quad
    {\bf p}'_0 := \left( {\bf I}- {\tt b}  {\tt b}^T  \diag({\bf Q})\right) {\bf p}'_{ref}(0),
\]
where ${\bf b}_i = 1 \ \forall i$. This choice for ${\bf p}^\prime(0)$ makes sure that ${\bf M}^\prime(0) = {\tt b}^T \diag({\bf Q}) {\tt p}^\prime(0) = 0$.

The relative root-mean-square error over time for the approximation on the uniform and curvilinear mesh are given in Table \ref{tab:uniform}. For each choice of interpolation spline, the best choice of mesh is emphasized in bold. As expected, the coarse interpolation results in the most accurate approximations if the mesh is coarse, and the fine interpolation when the grid is fine. The medium interpolation is most accurate when the grid is between coarse and fine. 

The errors increase slightly when time grows. However, the error percentages are still small enough to trust the approximation. Though the errors decrease rapidly with refinement of the mesh, the results do not show a clear rate of convergence. This is typical for the type of interpolation used.

Table \ref{tab:losses} shows the loss of mass and energy during the simulation with the medium interpolation method on a $20 \times 20$ uniform and curvilinear mesh. Though the accuracy is limited in this case (errors up to $4.47\%$ for the uniform mesh), the mass and energy losses are negligible up to machine accuracy. This is due to the symmetry-preserving nature of the discretization, and is true for all the simulations.

The results show that the method succeeds in constructing a symmetry-preserving discretization on a curvilinear, structured mesh.

\begin{table}[ht!]
\centering
\caption{Relative root-mean-square error (in percentages) of the approximation for different interpolation splines (taken from Table \ref{tab:splines}) and mesh sizes. For each choice of interpolation spline, the best choice of mesh is emphasized in bold. } \label{tab:uniform}

\vspace{0.2cm}

\resizebox{\textwidth}{!}{
\begin{tabular}{c|ccc|ccc|ccc}
\multicolumn{10}{c}{Uniform mesh on $[0,1] \times [0,1]$} \\
Time & \multicolumn{3}{c|}{Coarse interpolation}       &      \multicolumn{3}{|c|}{Medium interpolation}        &     \multicolumn{3}{c}{Fine interpolation} \\
\hline
  &           20x20      & 40x40 & 80x80 &      20x20 &        40x40     &   80x80 &   20x20&   40x40 &      80x80       \\  
  \hline
1 &        \textbf{0.30\% } & 0.056\% & 0.023\% &      0.64\%  &   \textbf{0.011\%} &  0.0042\% &   0.85\% &   0.018\% & \textbf{0.00045\%} \\
2 &        \textbf{0.59\% } & 0.111\% & 0.045\% &      1.24\%  &   \textbf{0.020\%} &  0.0055\% &   1.66\% &   0.032\% & \textbf{0.00057\%} \\
3 &        \textbf{0.86\% } & 0.171\% & 0.070\% &      1.81\%  &   \textbf{0.025\%} &  0.0053\% &   2.41\% &   0.052\% & \textbf{0.00061\%} \\
4 &        \textbf{1.11\% } & 0.223\% & 0.091\% &      2.32\%  &   \textbf{0.036\%} &  0.0070\% &   3.11\% &   0.063\% & \textbf{0.00057\%} \\
5 &        \textbf{1.33\% } & 0.276\% & 0.113\% &      2.75\%  &   \textbf{0.042\%} &  0.0066\% &   3.70\% &   0.081\% & \textbf{0.00051\%} \\
6 &        \textbf{1.53\% } & 0.335\% & 0.137\% &      3.14\%  &   \textbf{0.050\%} &  0.0071\% &   4.27\% &   0.096\% & \textbf{0.00041\%} \\
7 &        \textbf{1.71\% } & 0.390\% & 0.160\% &      3.46\%  &   \textbf{0.059\%} &  0.0088\% &   4.76\% &   0.113\% & \textbf{0.00058\%} \\
8 &        \textbf{1.84\% } & 0.447\% & 0.183\% &      3.80\%  &   \textbf{0.067\%} &  0.0104\% &   5.32\% &   0.127\% & \textbf{0.00060\%} \\
9 &        \textbf{1.94\% } & 0.507\% & 0.208\% &      4.11\%  &   \textbf{0.073\%} &  0.0101\% &   5.86\% &   0.146\% & \textbf{0.00070\%} \\
10 &       \textbf{2.00\% } & 0.557\% & 0.228\% &      4.47\%  &   \textbf{0.085\%} &  0.0124\% &   6.48\% &   0.160\% & \textbf{0.00064\%} \\
\hline
\multicolumn{10}{c}{Curvilinear mesh on $[0,1] \times [0,1]$} \\
Time & \multicolumn{3}{c|}{Coarse interpolation}       &      \multicolumn{3}{|c|}{Medium interpolation}        &     \multicolumn{3}{c}{Fine interpolation} \\
\hline
   &             20x20  &   40x40 &  80x80 &      20x20&      40x40      &   80x80 &    20x20 &  40x40 &    80x80         \\
   \hline
 1 &       \textbf{2.3\%} &   0.052\% &  0.024\% &       2.9\% &  \textbf{0.052\%} &  0.0045\% &      3.2\% &  0.079\% &  \textbf{0.0044\%} \\ 
 2 &       \textbf{3.7\%} &   0.098\% &  0.048\% &       4.8\% &  \textbf{0.093\%} &  0.0062\% &      5.4\% &  0.149\% &  \textbf{0.0053\%} \\
 3 &       \textbf{4.9\%} &   0.139\% &  0.070\% &       6.5\% &  \textbf{0.119\%} &  0.0071\% &      7.4\% &  0.207\% &  \textbf{0.0055\%} \\
 4 &       \textbf{6.9\%} &   0.191\% &  0.097\% &       9.1\% &  \textbf{0.168\%} &  0.0105\% &     10.1\% &  0.286\% &  \textbf{0.0081\%} \\
 5 &       \textbf{7.5\%} &   0.228\% &  0.119\% &      10.0\% &  \textbf{0.203\%} &  0.0119\% &     11.1\% &  0.354\% &  \textbf{0.0085\%} \\
 6 &       \textbf{9.4\%} &   0.282\% &  0.145\% &      12.2\% &  \textbf{0.251\%} &  0.0147\% &     13.5\% &  0.432\% &  \textbf{0.0108\%} \\
 7 &       \textbf{0.2\%} &   0.328\% &  0.169\% &      12.7\% &  \textbf{0.289\%} &  0.0165\% &     13.7\% &  0.497\% &  \textbf{0.0118\%} \\
 8 &       \textbf{0.6\%} &   0.366\% &  0.188\% &      13.0\% &  \textbf{0.315\%} &  0.0179\% &     13.9\% &  0.559\% &  \textbf{0.0120\%} \\
 9 &       \textbf{2.6\%} &   0.425\% &  0.219\% &      14.6\% &  \textbf{0.373\%} &  0.0210\% &     15.4\% &  0.645\% &  \textbf{0.0146\%} \\
 10 &      \textbf{1.8\%} &   0.455\% &  0.237\% &      12.7\% &  \textbf{0.388\%} &  0.0225\% &     13.0\% &  0.686\% &  \textbf{0.0147\%} \\

\end{tabular}
}
\end{table}

\begin{table}[ht!]
\centering
\caption{Percentage of mass and energy loss in the approximation on a $20 \times 20$ (coarse) mesh, using an interpolation spline with parameters {\tt nSpan} $ = 3$, {\tt nCont} $=1$, {\tt Order } $=11$, {\tt nConsist} $= 3$, {\tt wmax } $=0.6$ (medium interpolation). } \label{tab:losses}

\vspace{0.2cm}
\begin{tabular}{c|cc|cc}
& \multicolumn{2}{c|}{Uniform mesh} & \multicolumn{2}{c}{Curvilinear mesh} \\ 
\hline
Time & Mass loss & Energy loss & Mass loss & Energy loss \\
\hline
1 & 3.5E-08\% &  1.8E-11\% &  3.9E-07\% & 1.5E-12\% \\
2 & 7.0E-08\% &  2.5E-11\% &  8.3E-07\% & 2.2E-12\% \\
3 & 1.0E-07\% &  3.2E-11\% &  1.1E-06\% & 2.9E-12\% \\
4 & 1.4E-07\% &  3.8E-11\% &  1.6E-06\% & 3.8E-12\% \\
5 & 1.7E-07\% &  4.2E-11\% &  1.9E-06\% & 4.9E-12\% \\
6 & 2.1E-07\% &  4.4E-11\% &  2.4E-06\% & 6.0E-12\% \\
7 & 2.4E-07\% &  4.8E-11\% &  2.8E-06\% & 7.0E-12\% \\
8 & 2.8E-07\% &  5.2E-11\% &  3.1E-06\% & 8.1E-12\% \\
9 & 3.1E-07\% &  5.6E-11\% &  3.6E-06\% & 9.2E-12\% \\
10& 3.4E-07\% &  6.1E-11\% &  3.9E-06\% & 1.0E-11\% \\
\end{tabular}
\end{table}

\section{Conclusion}\label{sec:conclusion}
This paper describes a simple and effective strategy for the construction of symmetry-preserving 
methods on curvilinear, structured grids, offering flexibility and accuracy of the numerical approximations. The key idea is to use mutually-adjoint interpolation and sampling operators to switch from the continuous to the discrete operator. 
The numerical example shows that the method leads to results in which a high accuracy can be obtained, while the discrete mass and energy are 
both preserved.

The simulation of flow phenomena typically uses staggered grids, in which scalar fields and vector field 
components each have their own grids, and where each grid is shifted half a grid size with respect to the other grids.
Applying the discretization strategy explained in this paper is possible on such grids, but the details are outside the 
scope of this paper and will be the subject of a next paper.

One area of concern might be the number of non-zero elements in the discretization matrix. In a 3D calculation, any interpolation that is nonzero in a certain block results in a nonzero in the discretization matrix. Therefore, the number of nonzeros 
on every row will be ${\tt (4*nSpan-1)^3}$. The reduction of this number of nonzeros for each row will be discussed in another paper.

Finally, future work includes handling local grid refinements and the symmetry-preserving treatment of the compressible flow equations.

\bigskip
\textbf{Acknowledgments}

The authors gratefully wish to acknowledge the useful comments provided by Marc Gerritsma and Artur Palha from TU Delft and Eric Duveneck from Shell that helped to shape this work. We also thank VORtech and Shell that allowed us to develop this work.

\bibliographystyle{plain}
\bibliography{References}
\end{document}